\documentclass[12pt, leqno]{amsart}
\usepackage{geometry}
\usepackage[utf8]{inputenc}
\usepackage[unicode]{hyperref}
\usepackage{amsmath, amsthm, amscd, amssymb, amsfonts}
\usepackage{xcolor}
\colorlet{BLUE}{blue}
\usepackage{cancel}
\allowdisplaybreaks
\newcommand{\R}{\mathbb{R}}

\newcommand{\N}{\mathbb{N}}

\newcommand{\innerproduct}[2]{\left\langle #1, #2 \right\rangle}

\newtheorem{theorem}{Theorem}[section]
\newtheorem{lemma}[theorem]{Lemma}

\newtheorem{conjecture}[theorem]{Conjecture}

\theoremstyle{definition}
\newtheorem{definition}[theorem]{Definition}

\newtheorem{remark}[theorem]{Remark}

\title[Super-$p$-biharmonic and $p$-biharmonic submanifolds]{Liouville type results for super-$p$-biharmonic immersions and $p$-biharmonic submanifolds}

\author[T. D. Dinh, N. T. Dung, G. Wheeler]{Tien Dat Dinh, Nguyen Thac Dung, and Glen Wheeler}

\date{}

\address[Tien Dat Dinh]{Institute for Mathematics and its Applications, School of Mathematics and Applied Statistics, University of Wollongong, Northfields Ave, Wollongong, NSW 2500, Australia}
\email{tiendatd@uow.edu.au}

\address[Nguyen Thac Dung]{Department of Mathematics - Mechanics - Informatics, Vietnam National University, University of Sciences, Hanoi, No. 334, Nguyen Trai Road, Thanh Xuan, Hanoi, Vietnam}
\email{dungmath@gmail.com or dungmath@vnu.edu.vn}

\address[Glen Wheeler]{Institute for Mathematics and its Applications, School of Mathematics and Applied Statistics, University of Wollongong, Northfields Ave, Wollongong, NSW 2500, Australia}
\email{glenw@uow.edu.au}


\begin{document}
\begin{abstract}

    We introduce a weak notion of $\varepsilon$-super-$p$-biharmonic immersion into Euclidean space, formulated through a divergence inequality for the mean curvature vector. This class naturally extends the $\varepsilon$-superbiharmonic immersions previously introduced by the third author while remaining well-defined for all $p>1$, including the singular range $1<p<2$. We prove a Liouville-type rigidity theorem showing that every proper $\varepsilon$-super-$p$-biharmonic immersion is minimal provided
    \[
        \liminf_{r\to\infty} \frac1{r^{2}} \int_{f^{-1}(B_r)} |\vec{H}| ^{q}\,d\mu=0
    \]
    for some $q>p-\varepsilon$. This extends the corresponding superbiharmonic rigidity theorem to the $p$-biharmonic setting under a natural growth hypothesis.

    We also establish a minimality theorem for complete $p$-biharmonic immersions into Riemannian manifolds of non-positive sectional curvature under a polynomial growth assumption on the mean curvature. More precisely, if
    \[
        \int_{B_r(x_0)} | \vec{H} | ^{\,q}\,d\mu = O(r^d)
    \]
    for some $q\ge0$ and $d>0$, then the immersion is minimal. For $q=0$ this recovers the polynomial volume-growth theorem of Cao--Luo \cite{CL16}, while for $q>0$ it provides a polynomial-growth analogue of Han's $L^q$-integrability rigidity theorem \cite{H15}.

\end{abstract}
    
\maketitle

\section{Introduction}

In 1991, Chen \cite{C91} famously proposed the following conjecture.
\begin{conjecture}[Chen's conjecture] \label{Chen}
    Suppose $f:M^m\to\R^{m+n}$ is a biharmonic immersion, that is,
    \[ \Delta \vec{H} \equiv 0, \]
    then $\vec{H} \equiv 0$ and $f$ is minimal.
\end{conjecture}
Other generalizations of this conjecture have also been considered, namely
\begin{conjecture}[Generalized Chen's conjecture] \label{generalized Chen}
    Let $N^{m+n}$ be a Riemannian manifold with nonpositive sectional curvature and $f:M^m\to N^{m+n}$ be a biharmonic immersion, i.e.,
    \[
        \Delta \vec{H} + \sum_{i=1}^m R^N\bigl(\vec{H},df(e_i)\bigr)df(e_i) \equiv 0,
    \]
    then $\vec{H} \equiv 0$.
\end{conjecture}
Chen's and generalized Chen's conjectures have been extensively studied in the past few decades with many notable partial results. In particular, Chen \cite{C98} and Jiang \cite{J87} confirmed Conjecture \ref{Chen} for surfaces in $\R^3$ while Dimitri\'{c} \cite{D92} verified it for curves in $\R^n$. He also proved that Conjecture \ref{Chen} holds for pseudo-umbilic biharmonic immersions with $m \neq 4$ and $n$ arbitrary, as well as for hypersurfaces in Euclidean space if the immersion possesses at most two distinct principal curvatures. Very recently, by developing the technique introduced by Fu, Hong, and Zhan \cite{FHZ21}, which helped resolve Conjecture \ref{Chen} for hypersurfaces in $\R^5$, Fu, Hong, and Yang \cite{FHY26} proved an even more general result: Every biharmonic hypersurface in space forms of nonpositive sectional curvature is minimal. This settles Chen's Conjecture for hypersurfaces of Euclidean spaces. The high codimension case is still open in general. A recent work by Maeta \cite{M26} showed that rotationally symmetric biharmonic surfaces in $\R^4$ are minimal. For a more comprehensive survey on Chen's conjecture and related topics, readers are referred to a recent survey by Chen \cite{C14} and the references therein.

Surprisingly, Ou and Tang \cite{OT12} showed that Conjecture \ref{generalized Chen} does not hold in general. Their seminal paper suggests another active research direction, which is to find sufficient conditions for the minimality of biharmonic submanifolds. Luo \cite{L15} proved that the generalized Chen's conjecture is true for biharmonic submanifolds whose mean curvature vector is $L^q$ for some $q \in (0,+\infty)$ and the ambient space has non-positive sectional curvature. He further observed that, in the case of biharmonic hypersurfaces, the assumption that the ambient space has non-positive sectional curvature can be replaced by the weaker assumption that its Ricci curvature is non-positive. Other active lines of research include studying the $p$-biharmonic version of Chen's conjecture or considering the class of superbiharmonic submanifolds. The first question has been addressed by many authors, say, Han \cite{H15}, Han and Zhang \cite{HZ15}, and was confirmed for proper immersions by Cao and Luo\cite{CL16}. The second question was studied by the third author \cite{W13}. Recently, in \cite{BDST24}, Baker \textit{et al.} established rigidity results for $p$-biharmonic maps and $p$-biharmonic submanifolds, extending the line of investigation surrounding Chen's conjecture and its $p$-biharmonic generalizations. In particular, for $p$-biharmonic submanifolds, the authors gave sufficient conditions under which $p$-biharmonicity forces minimality, including smallness conditions on the $L^{n/2}$-norm of the sectional curvature and lower bounds for the fundamental tone. For $p=2$, these results recover rigidity results for biharmonic submanifolds and provide further partial affirmative results toward the generalized Chen conjecture.

Throughout this paper, we denote by $\vec{H}$ and $u := | \vec{H} |$ the mean curvature vector and its length, respectively.

\begin{theorem}[{\cite{W13}, Theorem 1.1}] \label{glen}
    Let $f: M^m \to N^{m+n}$ be an $\varepsilon$-superbiharmonic proper smooth immersion, i.e.,
    \begin{equation} \label{superbiharmonic}
        \innerproduct{\Delta \vec{H}}{\vec{H}} \geq (\varepsilon - 1) |\nabla\vec{H}|^2,
    \end{equation}
    for some $\varepsilon \in (0,1]$ where $\vec{H}$ is the mean curvature vector. Assume in addition that $f$ satisfies the growth condition
    \begin{equation}
        \liminf_{\rho \to \infty} \frac{1}{\rho^2} \int_{f^{-1}(B_\rho)} u^2 \,d\mu = 0,
    \end{equation}
    for every geodesic ball $B_\rho$ in $N^{m+n}$. Then $\vec{H} \equiv 0$ and $f$ is minimal.
\end{theorem}
In 2015, Han \cite{H15} showed that: Let $f: M^m \to N^{m+n}$ be a complete $p$-biharmonic immersion $( p \geq 2 )$ into a Riemannian manifold $N$ whose sectional curvature is nonpositive. If $u \in L^q(M)$ for some $q>0$, then $\vec{H}\equiv0$ and $f$ is minimal. He also proved the following result for a larger class of $\varepsilon$-super-$p$-biharmonic submanifolds. 
\begin{theorem}[{\cite{H15}, Theorem 3.4}] \label{Han's theorem}
    Let $f:M^m \to N^{m+n}$ be a complete immersion satisfying
    \begin{equation} \label{han's p-superbiharmonic}
        \innerproduct{\Delta (u^{p-2} \vec{H})}{ u^{p-2} \vec{H}} \geq (\varepsilon - 1) | \nabla ( u^{p-2} \vec{H} ) |^2
    \end{equation}
    for some $\varepsilon \in (0,1]$ and $p \geq 2$. If there is a constant $a \geq 0$ such that
    \[
        \int_M u^{a + 2p - 2} d\mu < +\infty,
    \]
    then $\vec{H}\equiv 0$ and $f$ is minimal.
\end{theorem}

In the present paper, we introduce the following notion of $\varepsilon$-super-$p$-biharmonic immersions into Euclidean space that is slightly different from \eqref{han's p-superbiharmonic}:

\begin{equation} \label{p-superbiharmonic}
    \innerproduct{\Delta(u^{p-2}\vec{H})}{\vec{H}} \geq (\varepsilon - 1) u^{p-2} \left( |\nabla^f \vec{H}|^2 + (p-2)|\nabla^M u|^2 \right).
\end{equation}

The detailed explanation for this choice of $\varepsilon$-super-$p$-biharmonic immersion and the comparison between this and Han's formulation is provided in Section 2. An important technical subtlety is that for $1 < p < 2$, the function $u^{p - 2}$ is singular wherever $\vec{H} = 0$. Hence, inequality \eqref{p-superbiharmonic} cannot be assumed to hold point-wise on $M$. A natural way to avoid this is to assume \eqref{p-superbiharmonic} only on $\Omega := \{ u > 0 \}$. However, to cover all variational $L^p$-biharmonic immersions with $1 < p < 2$, we use the following formulation.
    
\begin{definition} \label{weak formulation}
    Let $p > 1$ and $\varepsilon \in (0,1]$. A smooth immersion $f: M^m \to \R^{m + n}$ is termed $\varepsilon$-super-$p$-biharmonic if the divergence-form inequality
    \begin{equation}\label{weak-div}
        -\int_\Omega \langle X,\nabla\phi\rangle\, d\mu \geq \varepsilon\int_\Omega \phi\,u^{p-2} \left(|\nabla^f \vec{H}|^2+(p-2)|\nabla^M u|^2\right)\, d\mu
    \end{equation}
    holds for every nonnegative $\phi\in C_c^1(\Omega)$. Here $X = \innerproduct{\nabla^f ( u^{p - 2} \vec{H} )}{ \vec{H}}$.
\end{definition}

\begin{remark}
    It can be seen that \eqref{p-superbiharmonic} implies \eqref{weak-div} wherever all the quantities in \eqref{p-superbiharmonic} all well-defined. Thus, the notion used in Definition \ref{weak formulation} is sometimes called the weak formulation of \eqref{p-superbiharmonic}. Another reason for this choice of terminology is that \eqref{weak-div} is the distributional version of
    \[
        \text{div }X \geq \varepsilon u^{p - 2} \left( | \nabla^f \vec{H} |^2 + (p - 2) | \nabla^M u |^2 \right),
    \]
    with $X$ defined as in Definition \ref{weak formulation}.
\end{remark}

\begin{theorem} \label{main}
    Let $f: M^m \to \R^{m+n}$ be an $\varepsilon$-super-$p$-biharmonic proper smooth immersion with $p > 1$ and $\varepsilon \in (0,1]$. Let $\vec{H}$ denote the mean curvature vector and let $\nabla^f$ denote the pullback Euclidean connection along $f$. If, for some $q > p - \varepsilon$,
    \begin{equation} \label{growth condition}
        \liminf_{r \to \infty} \frac{1}{r^2} \int_{f^{-1}(B_r (x_0))} u^q \,d\mu = 0,
    \end{equation}
    where $B_r(x_0)$ denotes the ball of radius $r$ centered at any point $x_0 \in \R^{m + n}$, then $\vec{H} = 0$, and hence $f$ is minimal.
\end{theorem}
\begin{remark}
    It is crucial that $\nabla^f$ is the pull-back Euclidean connection because Theorem \ref{main} does not hold for the normal connection $\nabla^\perp$. Indeed, round spheres in Euclidean space are not minimal despite having parallel normal mean curvature.
\end{remark}

\begin{remark}
    The assumption $\varepsilon > 0$ is necessary. Again, consider round spheres in Euclidean space. Their mean curvature vector is nowhere vanishing but they satisfy both the growth condition \eqref{growth condition} for every $q > 0$ and
    \[
        \innerproduct{\Delta(u^{p-2}\vec{H})}{\vec{H}} = -u^{p-2} \bigl( |\nabla^f \vec{H}|^2 + (p-2)|\nabla^M u|^2 \bigr).
    \]
\end{remark}

The paper is organised as follows. In Section 2 we fix our notations, recall the $L^p$-bienergy convention, and explain the relation between our notion of $\varepsilon$-super-$p$-biharmonicity and Han's formulation. In Section 3, we prove Theorem \ref{main}, treating the classical formulation \eqref{p-superbiharmonic} on $\{ u > 0 \}$ and the divergence-form weak formulation \eqref{weak-div} separately. In Section 4, we adapt Han's local estimates to obtain a polynomial-growth minimality result for complete $p$-biharmonic immersions into Riemannian manifolds having nonpositive sectional curvature.

\section{Preliminaries}

\subsection{Notation and conventions}

    Throughout this paper, our objects of study are smooth immersions $f : M^m \to N^{m + n}$ with arbitrary codimension. In Sections 1--3, the ambient space is $\R^{m + n}$ and $f$ is further assumed to be proper. In Section 4, we consider $p$-biharmonic complete immersions into Riemannian manifolds with nonpositive sectional curvature. In both cases, the immersion induces on $M$ the structure of an $m$-dimensional Riemannian manifold. In particular, on each local coordinate chart $x^1,...,x^m$ of $M$, the Riemannian metric $g$ and the induced volume form are given locally by
\[
    g_{ij} = \innerproduct{\partial_i f}{\partial_j f}, \qquad d\mu = \sqrt{\text{det} g_{ij}}d\mathcal{L}^m,
\]
where $d\mathcal{L}^m$ is the $m$-dimensional Lebesgue (or Hausdorff) measure. Associated with the immersion $f$ is the vector-valued second fundamental form, defined as
\[
    A_{ij} = (\overline{\nabla}_i\partial_j f)^\perp,
\]
where $\overline{\nabla}$ is the covariant derivative with respect to the Levi-Civita connection of $N$ and $(\cdot)^\perp$ is the projection onto the normal bundle $TM^\perp$,
\[
    X^\perp = X - X^\top = X - g^{ij}\innerproduct{X}{\partial_i f}\partial_j f
\]
for every vector field $X: M^m \to TN^{m+n}$. The mean curvature vector is then defined as the trace of $A$ under the metric $g$,
\[
    \vec{H} := g^{ij}A_{ij}.
\]
Recall the Gauss formula,
\[
    \overline{\nabla}_XY = \nabla^M_X Y + A(X , Y), 
\]
which says that the Levi-Civita connection $\nabla^M$ corresponding to the metric $g$ is the tangential projection of the ambient connection $\overline{\nabla}$ of $N$. In a local coordinate $x^1,...,x^m$ of $M$, the connection $\nabla^M$ is given by
$$
    \nabla^M_{\partial_i}\partial_j = \Gamma_{ij}^k\partial_k, \qquad \text{where } \Gamma_{ij}^k = \frac{1}{2} g^{kl} (\partial_i g_{jl} + \partial_j g_{il} - \partial_l g_{ij}).
$$

For a vector field $W\in\Gamma(f^*TN)$ along the immersion, we write
\[
    \nabla^f_XW:=\overline{\nabla}_{df(X)}W
\]
for the pull-back connection and reserve $\nabla^M$ for the intrinsic Levi--Civita connection on $M$. The corresponding connection Laplacian on $f^*TN$ is
\[
    \Delta W:=\sum_{i=1}^m\left(\nabla^f_{e_i}\nabla^f_{e_i}W-\nabla^f_{\nabla^M_{e_i}e_i}W\right)=-(\nabla^f)^*\nabla^fW.
\]
For scalar functions, $\nabla$ and $\Delta$ denote the intrinsic gradient and Laplace--Beltrami operator. For a normal field $\xi$, the Weingarten formula
\[
    \nabla^f_X\xi=-df(S_\xi X)+\nabla^\perp_X\xi,
    \qquad
    \langle S_\xi X,Y\rangle=\langle A(X,Y),\xi\rangle,
\]
defines the shape operator $S_\xi$. When $N=\R^{m+n}$, regarding $f$ as a section of $f^*T\R^{m+n}$ gives $\Delta f=\vec H$.

\subsection{Preliminaries on \textit{p}-biharmonic maps}

The following variational approach allows us to generalize the above definition to $p$-biharmonicity (a comprehensive treatment of the subject can be found in \cite{C18}, \cite{LC26} and the references therein). Let $\varphi: M^m \to N^{m + n}$ be a smooth map between Riemannian manifolds and $p > 1$. We say that $\varphi$ is $p$-biharmonic if it is critical, with respect to every compactly supported smooth variation, for the $L^p$-bienergy functional
\begin{equation} \label{p-bienergy}
    E(\varphi) := \frac{1}{p} \int_D |\tau(\varphi)|^p \,d\mu,
\end{equation} 
where $D\Subset M$ is any relatively compact domain containing the support of the variation. Here
$$
\tau(\varphi) := \text{tr}_g \nabla d\varphi = \nabla^\varphi_{e_i} d\varphi(e_i) - d\varphi(\nabla_{e_i}^M e_i)
$$
is the tension field. The Euler-Lagrange equation of \eqref{p-bienergy} is given by $\tau_p(\varphi) = 0$ with $\tau_p(\varphi)$ being the $p$-bitension field of $\varphi$, defined as
\begin{equation} \label{phan}
    \tau_p(\varphi) := \Delta\bigl( |\tau(\varphi)|^{p-2}\tau(\varphi) \bigr) + \sum_{i=1}^m R^N(|\tau(\varphi)|^{p-2}\tau(\varphi), d\varphi(e_i)) d\varphi(e_i),
\end{equation}
where $R^N$ is the Riemann curvature tensor of $N^{m + n}$, and $\Delta$ is the pull-back connection Laplacian on $\varphi^*TN$. When considering immersions into Euclidean space, the ambient space is flat, $R^N\equiv0$, and hence the curvature term in $\tau_p(\phi)$ vanishes. Thus, the $p$-biharmonicity of such an immersion can be written as
\[
    \Delta_p\vec{H} := \Delta\bigl( u^{p-2}\vec{H} \bigr) = 0,
\]
where $\vec{H}$ is the mean curvature vector. In particular, the manifold $M^m$ is itself called a $p$-biharmonic submanifold of $\R^{m+n}$.

It is worth noting that we are using the $L^p$-bienergy convention, which is distinct from those based on $\tau_p$, termed bi-$p$-harmonic maps, which were explicitly defined and studied in \cite{LC26}. In fact, it is interesting to ask whether these results still hold true in that context.

\begin{remark}
    It appears that two different curvature conventions are used in \cite{CL16}.  In (2.1), the authors adopt the curvature convention used in  \cite{EL83}, whereas in the derivation of (2.6) and (2.7), the curvature tensor $R^N$ differs by a minus sign. With the latter convention, formulas (2.6) and (2.7) in \cite{CL16} coincide with formulas (2.2) and (2.3) in \cite{H15}, respectively. Throughout this paper, we use the curvature convention similar to that in \cite{J87, H15}.
\end{remark}

\subsection{Super-$p$-biharmonic immersions}

In the pursuit of extending the theory of $\varepsilon$-super biharmonic maps, there appears to be many different choices for the suitable notion of $\varepsilon$-super-$p$-biharmonicity. When $p = 2$, both \eqref{han's p-superbiharmonic} and \eqref{p-superbiharmonic} reduce to the standard $\varepsilon$-superbiharmonic condition \eqref{superbiharmonic} established in \cite{W13}. However, since these formulations differ for $p \neq 2$, our definition describes a class distinct from the class considered by Han \cite{H15}. In this subsection, we would like to explain our choice for \eqref{p-superbiharmonic} and \eqref{weak-div}. In the following computation, we abbreviate the pull-back connection $\nabla^f$ by $\nabla$.

Our formulation based on the standard Bochner type identity for the connection Laplacian of the squared norm, i.e.,
\[
    \frac{1}{2}\Delta u^2 = \innerproduct{\Delta \vec{H}}{\vec{H}} + |\nabla \vec{H} |^2,
\]
where $u := | \vec{H} |$. Substituting \eqref{superbiharmonic} into the above identity yields $ \frac{1}{2}\Delta u^2 \geq \varepsilon |\nabla \vec{H} |^2$, which eventually implies $\Delta u ^2 \geq 0$ when $\varepsilon = 0$. For our purpose, we need to find the corresponding Bochner type identity for $u^p$. Consider an orthonormal frame $e_1,...,e_m$ of $M$. The product rule of covariant differentiation yields
\[
    \nabla_i \innerproduct{\nabla_i (u^{p-2}\vec{H})}{\vec{H}} = \innerproduct{\nabla_i \nabla_i (u^{p-2}\vec{H})}{\vec{H}} + \innerproduct{\nabla_i (u^{p-2}\vec{H})}{\nabla_i \vec{H}}.
\]
Since
\[
    \nabla_i (u^{p-2}\vec{H}) = \nabla_i (u^{p-2}) \vec{H} + u^{p-2} \nabla_i \vec{H},
\]
we have that
\[
    \begin{aligned}
        \innerproduct{\nabla_i (u^{p-2}\vec{H})}{\vec{H}}
        &= \nabla_i (u^{p-2}) u^2 + u^{p-2} \innerproduct{\nabla_i \vec{H}}{\vec{H}} \\
        &= (p-1) u^{p-2} \innerproduct{\nabla_i \vec{H}}{\vec{H}}. 
    \end{aligned}  
\]
Hence,
\[
    \begin{aligned}
        \nabla_i \innerproduct{\nabla_i (u^{p-2}\vec{H})}{\vec{H}} &= \nabla_i \left( (p-1)u^{p-2}\langle \vec{H}, \nabla_i \vec{H} \rangle \right) \\
        &= \frac{p-1}{p} \Delta u^p.
    \end{aligned}
\]
On the other hand,
\[
    \begin{aligned}
        \innerproduct{ \nabla_i u^{p-2} \vec{H}}{\nabla_i \vec{H}}
        &= \innerproduct{(p-2) u^{p-4} \innerproduct{\vec{H}}{\nabla_i \vec{H}} \vec{H} + u^{p-2} \nabla_i \vec{H}}{\nabla_i \vec{H}}\\
        &= (p-2) u^{p-4} \innerproduct{\vec{H}}{\nabla_i \vec{H}}^2 + u^{p-2} |\nabla \vec{H} |^2\\
        &= u^{p-2} \left( |\nabla \vec{H}|^2 + (p-2)|\nabla u|^2 \right),
    \end{aligned}
\]
where the third equality is due to $\sum_{i=1}^m\innerproduct{\vec{H}}{\nabla_i \vec{H}}^2 = u^2 |\nabla u|^2$. Finally, we equate corresponding terms and obtain
\[
    \frac{p-1}{p} \Delta u^p = \innerproduct{\Delta(u^{p-2}\vec{H})}{\vec{H}} + u^{p-2} \left( |\nabla \vec{H}|^2 + (p-2)|\nabla u|^2 \right). 
\]
With respect to the above identity, we introduce  inequality \eqref{p-superbiharmonic} based on the idea of the case $p = 2$.

\begin{remark}
    For $1 < p < 2$, the vector field $u^{p - 2} \vec{H}$ is singular at points where $\vec{H}$ vanishes. Thus a smooth immersion does not automatically give a classical equation for $\Delta(u^{p-2}\vec{H})$ across the zero set. Our first formulation \eqref{p-superbiharmonic} avoids this by assuming the inequality only on $\{\vec{H}\ne0\}$. However, in order to cover all variational $L^p$-biharmonic immersions with $1 < p < 2$, we use the weak formulation \eqref{weak-div}, which is the distributional version of
    \[
        \text{div} X \geq \varepsilon u^{p - 2} \left( | \nabla \vec{H} |^2 + (p - 2) | \nabla u |^2 \right),
    \]
    where $X := \innerproduct{ \nabla (u^{p-2} \vec{H})}{\vec{H}}$.
\end{remark}

\begin{remark}
    By direct computation, one can show that on $\{ u > 0 \}$,
    \[
        \left|\nabla\big(u^{p-2}\vec{H}\big)\right|^2 = u^{2p-4} \left( |\nabla\vec{H}|^2+p(p-2)|\nabla u|^2 \right).
    \]
    Therefore,
    \[
        u^{2-p} \left|\nabla\big(u^{p-2}\vec{H}\big)\right|^2 = u^{p-2} \left( |\nabla\vec{H}|^2+p(p-2)|\nabla u|^2 \right).
    \]
    After dividing Han's condition \eqref{han's p-superbiharmonic} by $u^{p-2}$, the comparison between his condition with our proposed condition \eqref{p-superbiharmonic} is governed by the corresponding coefficients
    \[
        p(p - 2) \qquad \text{and} \qquad p-2.
    \]
    Thus, for $p \geq 2$, our condition is stronger than Han's condition. For $ 1 < p < 2$, the direction is reversed.
\end{remark}
\section{A Liouville type result for $\varepsilon$-super-$p$-biharmonic immersions into Euclidean space}

In this section, we use the abbreviated notation $\nabla$ for the pull-back connection.

\begin{lemma} \label{firstlemma}
Let $f:M^m\to\R^{m+n}$ be a smooth immersion. Denote by $\vec{H}$ the corresponding mean curvature vector of $M$ and $\nabla$ the pullback Euclidean connection along $f$. Then
    \[
        \innerproduct{\nabla(u^{p-2}\vec{H})}{\vec{H}} = (p-1)u^{p-1}\nabla u
    \]
    wherever $u > 0$.
\end{lemma}

\begin{proof}
    The product rule of covariant derivative shows that
    \[
        \nabla(u^{p-2}\vec{H}) = \nabla(u^{p-2}) \otimes \vec{H} + u^{p-2}\nabla\vec{H} = (p-2)u^{p-3}\nabla u \otimes \vec{H} + u^{p-2}\nabla\vec{H}.
    \]
    Thus,
    \[
        \innerproduct{\nabla(u^{p-2}\vec{H})}{\vec{H}} = (p-2)u^{p-3} \innerproduct{\nabla u \otimes \vec{H}}{\vec{H}} + u^{p-2}\innerproduct{\nabla\vec{H}}{\vec{H}}.
    \]
    Since $u^2 = \innerproduct{\vec{H}}{\vec{H}}$, we have that $u\nabla u = \innerproduct{\vec{H}}{\nabla\vec{H}}$. Therefore,
    \[
        \innerproduct{\nabla(u^{p-2}\vec{H})}{\vec{H}} = (p-2)u^{p-3}(u^2 \nabla u) + u^{p-2}(u\nabla u) = (p-1)u^{p-1}\nabla u.
\qedhere    \]
\end{proof}

Next, we calculate the divergence of the vector field defined in Lemma \ref{firstlemma}.

\begin{lemma} \label{secondlemma}
Let $f:M^m\to\R^{m+n}$ be a smooth immersion. Denote by $\vec{H}$ the corresponding mean curvature vector of $M$ and $\nabla$ the pullback Euclidean connection along $f$. Then
    \[
        \text{div} \innerproduct{\nabla(u^{p-2}\vec{H})}{\vec{H}} = \innerproduct{\Delta_p\vec{H}}{\vec{H}} + (p-2)u^{p-2}|\nabla u|^2 + u^{p-2}|\nabla\vec{H}|^2
    \]
    wherever $u > 0$.
\end{lemma}

\begin{proof}
    By the product rule of the divergence operator, we have
    \[
        \begin{aligned}
            \text{div} \innerproduct{\nabla(u^{p-2}\vec{H})}{\vec{H}} 
                &= \innerproduct{\text{tr}\nabla^2(u^{p-2}\vec{H})}{\vec{H}} + \innerproduct{\nabla(u^{p-2}\vec{H})}{\nabla\vec{H}} \\
                &= \innerproduct{\Delta_p\vec{H}}{\vec{H}} + \innerproduct{\nabla(u^{p-2}\vec{H})}{\nabla\vec{H}}.
        \end{aligned}
    \]
    As in the proof of Lemma \ref{firstlemma}, we use the product rule of covariant differentiation and obtain
    \[
        \innerproduct{\nabla(u^{p-2}\vec{H})}{\nabla\vec{H}} = (p-2)u^{p-3}\innerproduct{\nabla u \otimes \vec{H}}{\nabla\vec{H}} + u^{p-2}|\nabla\vec{H}|^2.
    \]
    Since
    \[
        \begin{aligned}
            \innerproduct{\nabla u \otimes \vec{H}}{\nabla\vec{H}} 
                &= \sum_{i=1}^m \nabla_iu \innerproduct{\vec{H}}{\nabla_i\vec{H}} \\
                &= \sum_{i=1}^m \nabla_iu \bigl( u\nabla_iu \bigr) \\
                &= u|\nabla u|^2,
        \end{aligned}
    \]
    we finally have
    \[
        \text{div} \innerproduct{\nabla(u^{p-2}\vec{H})}{\vec{H}} = \innerproduct{\Delta_p \vec{H}}{\vec{H}} + (p-2)u^{p-2}|\nabla u|^2 + u^{p-2}|\nabla\vec{H}|^2.
\qedhere    \]
\end{proof}

In the sequel, we will employ cut-off functions defined as follows: Let $\tilde{\eta}$ be a smooth cut-off function on a ball $B_{r}(x_0)$ in $\R^{m+n}$ centered at any point $x_0 \in \R^{m + n}$ such that $\tilde{\eta} \equiv 1$ in $B_r(x_0)$, $\tilde{\eta} \equiv 0$ outside $B_{2r}(x_0)$, and $|\nabla\tilde{\eta}| < \frac{c}{r}$ for some positive constant $c$. Define
\begin{equation} \label{cut-off function}
    \eta = \tilde{\eta} \circ f.
\end{equation}
Since $f$ is proper, $\eta \in C_c^\infty(M)$ and since $f$ is an isometric immersion for the induced metric, $|\nabla^M\eta| < \frac{c}{r}$.

\begin{lemma} \label{main lemma}
    Let $f: M^m \to \R^{m+n}$ be a proper smooth immersion, $p > 1$ and $\varepsilon \in (0,1]$. Let $\vec{H}$ denote the mean curvature vector and let $\nabla$ denote the pullback Euclidean connection along $f$. Suppose that \eqref{p-superbiharmonic} holds on the set $\Omega := \{ u > 0 \}$ in the classical sense. If, for some $q > p - \varepsilon$,
    \[
        \liminf_{r \to \infty} \frac{1}{r^2} \int_{f^{-1}(B_r(x_0))} u^q \,d\mu = 0,
    \]
    where $B_r(x_0)$ is a ball or radius $r$, centered at any point $x_0 \in \R^{m + n}$, then $\nabla\vec{H} = 0$.
\end{lemma}
\begin{proof}[Proof of Lemma \ref{main lemma}]
    
    To handle the potential singularity across the zero set of $u$, consider super-level sets
    \[
        \Omega_\sigma := \{ u > \sigma \} \qquad \text{for regular values } \sigma > 0.
    \]
     Note that almost every $\sigma > 0$ is a regular value of $u$ so by Sard's Theorem, $\partial\Omega_\sigma$ is a smooth submanifold of $M$ for every such $\sigma$. In what follows, one first integrates over $\Omega_\sigma\cap K$, where $K$ is a compact set containing the support of $\eta$, and then observes that no extra boundary term appears away from $\partial\Omega_\sigma$.
     
     Let $V := u^{p - 2} \vec{H}$ and $\eta$ be as in \eqref{cut-off function}. Integrating $\text{div} (\eta^2 u^{q - p} \innerproduct{\nabla V}{\vec{H}})$ over $\Omega_\sigma$ and applying the Divergence Theorem, we obtain
    \[
        \begin{aligned}
            \int_{\Omega_\sigma} \text{div} (\eta^2 u^{q - p} \innerproduct{\nabla V}{\vec{H}}) \, d\mu &= \int_{\partial\Omega_\sigma} \eta^2 u^{q - p}\innerproduct{\nabla V}{\vec{H}} \cdot \nu \, dS \\
            &= - (p - 1) \sigma^{q - 1} \int_{\partial \Omega_\sigma} \eta^2 | \nabla u | \, dS \\
            & \leq 0,
        \end{aligned}
    \]
    where the second equality follows from Lemma \ref{firstlemma}. Here, $\nu := - \frac{\nabla u}{| \nabla u |}$ is the outward unit normal of $\partial \Omega_\sigma$. Now,
    \begin{align*}
        \int_{\Omega_\sigma}\text{div} (\eta^2 u^{q - p} \innerproduct{\nabla V}{\vec{H}}) \, d\mu  =& \int_{\Omega_\sigma}\left\langle \innerproduct{\vec{H}}{\nabla (u^{p-2}\vec{H})} , \nabla\bigl( \eta^2 u^{q-p} \bigr) \right\rangle \, d\mu \\
        &+ \int_{\Omega_\sigma} u^{q-p} \eta^2 \innerproduct{\vec{H}}{\Delta_p\vec{H}} \, d\mu\\
        &+ \int_{\Omega_\sigma} u^{q-p} \eta^2 \innerproduct{\nabla(u^{p-2}\vec{H})}{\nabla\vec{H}} \, d\mu \\
        =& (q - p)(p - 1) \int_{\Omega_\sigma}  u^{q-2} \eta^2 |\nabla u|^2 \, d\mu \\
        &+ 2(p - 1) \int_{\Omega_\sigma}  u^{q - 1} \eta \innerproduct{\nabla u}{\nabla\eta} \, d\mu \\
        &+ \int_{\Omega_\sigma} u^{q - p} \eta^2 \innerproduct{\vec{H}}{\Delta_p \vec{H}} \, d\mu + (p - 2) \int_{\Omega_\sigma}  u^{q - 2} \eta^2 | \nabla u |^2 \, d\mu \\
        &+ \int_{\Omega_\sigma} u^{q - 2} \eta^2 | \nabla \vec{H} |^2 \, d\mu,
    \end{align*}
    where the second equality follows from Lemma \ref{firstlemma} and Lemma \ref{secondlemma}. Using the $\varepsilon$-super $p$-biharmonicity of $f$, we have
    \[
        \begin{aligned}
            (\varepsilon - 1) \int_{\Omega_\sigma} u^{q-2}  \eta^2 \bigl( |\nabla\vec{H}|^2 +& (p-2)|\nabla u|^2 \bigr) \,d\mu \\
                \leq& \int_{\Omega_\sigma} u^{q-p} \eta^2 \innerproduct{\vec{H}}{\Delta_p\vec{H}} \,d\mu \\
                \leq& -(q-p)(p-1) \int_{\Omega_\sigma} u^{q-2} \eta^2 |\nabla u|^2 \,d\mu \\
                &- 2(p-1) \int_{\Omega_\sigma} u^{q-1} \eta \innerproduct{\nabla u}{\nabla\eta} \,d\mu \\
                &- (p-2) \int_{\Omega_\sigma} u^{q-2} \eta^2 |\nabla u|^2 \,d\mu \\
                &- \int_{\Omega_\sigma} u^{q-2} \eta^2 |\nabla\vec{H}|^2 \,d\mu.
        \end{aligned}
    \]
    Isolating coercive terms, we obtain
    \begin{equation} \label{first}
        \begin{aligned}
            &\varepsilon\int_{\Omega_\sigma} \eta^2 u^{q-2}|\nabla \vec H|^2\,d\mu + \varepsilon(p-2)\int_{\Omega_\sigma} \eta^2u^{q-2}|\nabla u|^2\,d\mu  \\
            &\qquad + (p-1)(q-p)\int_{\Omega_\sigma} \eta^2u^{q-2}|\nabla u|^2\,d\mu \\
            &\leq -2(p-1)\int_{\Omega_\sigma} \eta u^{q-1} \innerproduct{\nabla u}{\nabla \eta}\,d\mu.
        \end{aligned}
    \end{equation}
    Equivalently,
    \[
    \begin{aligned}
        \int_{\Omega_\sigma} \eta^2u^{q-2} \Big[ \varepsilon|\nabla \vec{H}|^2 + \big(\varepsilon(p-2)+(p-1)(q-p)\big)|\nabla u|^2 \Big]\,d\mu& \\
        \leq -2(p-1)\int_{\Omega_\sigma} \eta u^{q-1} \langle \nabla u,\nabla\eta\rangle\,d\mu.&
    \end{aligned}
    \]
    By the Peter-Paul inequality, the right-hand side satisfies the following estimate for every $\delta > 0$,
    \begin{equation} \label{second}
        \begin{aligned}
            -2(p-1) &\int_{\Omega_\sigma} u^{q-1} \eta \innerproduct{\nabla u}{\nabla\eta} \,d\mu \\
            &\leq \delta \int_{\Omega_\sigma} u^{q - 2} \eta^2 |\nabla u|^2 \,d\mu + C_\delta \int_{\Omega_\sigma} u^q |\nabla\eta|^2 \,d\mu,
        \end{aligned}
    \end{equation}
    with $C_\delta = \frac{(p - 1)^2}{\delta}$. Hence,
    \[
        \int_{\Omega_\sigma} \eta^2u^{q-2} \Big[ \varepsilon|\nabla \vec{H}|^2 + A_\delta |\nabla u|^2 \Big]\,d\mu \leq C_\delta\int_{\Omega_\sigma} u^q|\nabla\eta|^2\,d\mu ,
    \]
    where
    \[
        A_\delta:=\varepsilon(p-2)+(p-1)(q-p)-\delta .
    \]
    If \(A_\delta\geq0\), then the left-hand side directly controls
    \[
        \varepsilon\int_{\Omega_\sigma} | \nabla \vec{H} |^2 \, d\mu.
    \]
    If \(A_\delta<0\), Kato's inequality $|\nabla u|\leq |\nabla \vec{H}|$ gives
    \[
        \varepsilon|\nabla\vec H|^2+A_\delta|\nabla u|^2 \geq (\varepsilon+A_\delta)|\nabla\vec H|^2 = \big((p-1)(q-p+\varepsilon)-\delta\big)|\nabla\vec H|^2.
    \]
    Thus, provided $q > p - \varepsilon$, we may choose $0<\delta<(p-1)(q-p+\varepsilon)$ and obtain
    \begin{equation} \label{second 2}
        \gamma_\delta \int_{\Omega_\sigma} \eta^2u^{q-2}|\nabla\vec H|^2\,d\mu \leq C_\delta\int_{\Omega_\sigma} u^q|\nabla\eta|^2\,d\mu ,
    \end{equation}
    with
    \[
        \gamma_\delta := \min\left\{ \varepsilon,\, (p-1)(q-p+\varepsilon)-\delta \right\}>0 .
    \]
    We handle the left-hand side of \eqref{second 2} using Fatou's lemma:
    \[
        \gamma_\delta \int_\Omega u^{q-2} \eta^2 |\nabla\vec{H}|^2 \,d\mu \leq \liminf_{\sigma \to 0^+} \gamma_\delta \int_{\Omega_\sigma} u^{q-2} \eta^2 |\nabla\vec{H}|^2 \,d\mu.
    \]
    For the right-hand side, we argue as follows: Since $u^q$ is continuous on $\Omega$, the function $u^q | \nabla \eta |^2$ is bounded and integrable on the compact support of $\nabla \eta$. Moreover, because $\chi_{\Omega_\sigma}u^q |\nabla \eta|^2 \geq0$ and $\chi_{\Omega_\sigma}u^q|\nabla \eta|^2 \uparrow \chi_{\Omega}u^q|\nabla \eta|^2$ almost everywhere as $\sigma\to0^+$, the Monotone Convergence Theorem implies that
    \[
        \lim_{\sigma \to 0^+} \int_{\Omega_\sigma} u^q |\nabla\eta|^2 \,d\mu = \int_{\Omega} u^q |\nabla\eta|^2 \,d\mu.
    \]
    In short, letting $\sigma \to 0^+$ and then using the bound $|\nabla\eta| < \frac{c}{r}$ yields
    \begin{equation} \label{third}
        \begin{aligned}
            \int_\Omega u^{q-2} \eta^2 |\nabla\vec{H}|^2 \,d\mu 
                &\leq \frac{C(p,q,\varepsilon,\delta)}{\delta\gamma_\delta} \int_\Omega u^q |\nabla\eta|^2 \,d\mu \\
                &\leq \frac{C'(p,q,\varepsilon,\delta)}{r^2} \int_{f^{-1}(B_{2r}(x_0))} u^q \,d\mu.
        \end{aligned}
    \end{equation}
    To use the growth condition \eqref{growth condition}, consider a sequence $\{\rho_j\}_{j \in \N}$ of positive real numbers such that $\rho_j \to \infty$ and
    \[
        \lim_{j \to \infty} \frac{1}{\rho_j^2} \int_{f^{-1}(B_{\rho_j}(x_0))} u^q \,d\mu = 0.
    \]
    Taking $r_j = \rho_j/2$ in \eqref{third} gives us
    \[
        \lim_{j \to \infty} \int_{\Omega \cap f^{-1}(B_{r_j}(x_0))} u^{q-2} |\nabla\vec{H}|^2 \,d\mu = 0.
    \]
    After passing to an increasing sequence if necessary, these domains exhaust $\Omega$, which shows that
    \[
        \int_\Omega u^{q-2} |\nabla\vec{H}|^2 \,d\mu = 0.
    \]
    On $\Omega$, the weight $u^{q-2}$ is positive so $\nabla\vec{H} \equiv 0$ there. On the interior of $\{ u = 0 \}$, $\vec{H}$ is identically zero so $\nabla \vec{H} = 0$ there as well. At the boundary points of the zero set, the conclusion follows by the continuity of $\nabla \vec{H}$. Therefore, $\nabla\vec{H} \equiv 0$ on all of $M$.
\end{proof}

\begin{lemma} \label{main lemma, weak}
    Let $f: M^m \to \R^{m+n}$ be a proper smooth immersion, $p > 1$ and $\varepsilon \in (0,1]$. Let $\vec{H}$ denote the mean curvature vector and let $\nabla$ denote the pullback Euclidean connection along $f$. Suppose that \eqref{p-superbiharmonic} holds weakly on $\Omega := \{ u > 0 \}$, i.e., the divergence-form inequality \eqref{weak-div} holds for every nonnegative $\phi\in C_c^1(\Omega)$. If, for some $q > p - \varepsilon$,
    \[
        \liminf_{r \to \infty} \frac{1}{r^2} \int_{f^{-1}(B_r(x_0))} u^q \,d\mu = 0,
    \]
    where $B_r(x_0)$ is a ball or radius $r$, centered at any point $x_0 \in \R^{m + n}$, then $\nabla\vec{H} = 0$.
\end{lemma}

\begin{proof}
    For each regular value $\sigma > 0$ of $u$, let $\zeta_\sigma:[0,\infty)\to[0,1]$ be a smooth cut-off function  such that
    \[
        \zeta_\sigma(t)=0 \text{ for } t\leq\sigma, \qquad \zeta_\sigma(t)=1 \text{ for } t\geq 2\sigma, \qquad \zeta_\sigma'(t)\geq0.
    \]
    Such a function can be constructed as follows: choose one fixed nondecreasing smooth function
    \[
        \zeta:[0,\infty)\to[0,1], \qquad \zeta(t)=0 \text{ for } t\le1, \qquad \zeta(t)=1 \text{ for } t\ge2,
    \]
    and set
    \[
    \zeta_\sigma(t)=\zeta(t/\sigma).
    \]
    Then define
    \[
    \phi_\sigma=\eta^2u^{q-p}\zeta_\sigma(u),
    \]
    with $\eta$ be as in \eqref{cut-off function}. Each $\phi_\sigma$ is compactly supported in $\Omega$ because $\zeta_\sigma(u)$ vanishes near the zero set of $u$. Testing \eqref{weak-div} with $\phi_\sigma$ then yields:
    \begin{equation} \label{first1}
        -\int_\Omega \langle X,\nabla\phi_\sigma\rangle\, d\mu \geq \varepsilon\int_\Omega \phi_\sigma \,u^{p-2} \left(|\nabla\vec{H}|^2+(p-2)|\nabla u|^2\right)\, d\mu,
    \end{equation}
    with $X = \innerproduct{\nabla ( u^{p - 2} \vec{H} )}{ \vec{H}}$. Expanding the left-hand side of \eqref{first1}, we obtain
    \begin{align*}
        -\int_\Omega \langle X,\nabla\phi_\sigma\rangle\, d\mu &= -2(p-1)\int_\Omega \eta u^{q-1}\zeta_\sigma(u) \langle\nabla u,\nabla\eta\rangle\, d\mu \\
        &\quad -(p-1)(q-p)\int_\Omega \eta^2u^{q-2}\zeta_\sigma(u)|\nabla u|^2\, d\mu \\
        &\quad -(p-1)\int_\Omega \eta^2u^{q-1}\zeta_\sigma'(u)|\nabla u|^2\, d\mu .
    \end{align*}
    Because $\zeta_\sigma' \geq 0$, the last term on the right-hand side is nonpositive and it may be discarded. Substituting this into \eqref{first1} and isolating coercive terms, we obtain
    \begin{align*}
        &\int_\Omega \eta^2u^{q-2}\zeta_\sigma(u) \left[ \varepsilon|\nabla\vec{H}|^2+ \{\varepsilon(p-2)+(p-1)(q-p)\}|\nabla u|^2 \right]\, d\mu \\
        &\qquad\leq -2(p-1)\int_\Omega \eta u^{q-1}\zeta_\sigma(u) \langle\nabla u,\nabla\eta\rangle\, d\mu .
    \end{align*}
    We then proceed exactly as in the proof of Lemma \ref{main lemma}, that is, using Peter-Paul inequality and Kato inequality along side the condition $q > p - \varepsilon$ and get
    \[
        \gamma_\delta\int_\Omega \eta^2u^{q-2}\zeta_\sigma(u)|\nabla\vec{H}|^2\, d\mu \leq C_\delta\int_\Omega u^q\zeta_\sigma(u)|\nabla\eta|^2\, d\mu .
    \]
    Then $\zeta_\sigma(u)\uparrow1$ on $\Omega$ as \(\sigma\downarrow0 \) so by applying Fatou's lemma for the left-hand side and Monotone Convergence Theorem for the right, we have, by letting $\sigma \to 0^+$, that
    \begin{equation} \label{fifth}
        \begin{aligned}
            \int_\Omega u^{q-2} \eta^2 |\nabla\vec{H}|^2 \,d\mu 
                &\leq \frac{C(p,q,\varepsilon,\delta)}{\delta\gamma_\delta} \int_\Omega u^q |\nabla\eta|^2 \,d\mu \\
                &\leq \frac{C'(p,q,\varepsilon,\delta)}{r^2} \int_{f^{-1}( B_{2r}(x_0))} u^q \,d\mu.
        \end{aligned}
    \end{equation}
    To use the growth condition \eqref{growth condition}, let $\{\rho_j\}_{j \in \N}$ be a sequence of positive real numbers such that $\rho_j \to \infty$ and
    \[
        \lim_{j \to \infty} \frac{1}{\rho_j^2} \int_{f^{-1}(B_{\rho_j}(x_0))} u^q \,d\mu = 0.
    \]
    Taking $r_j = \rho_j/2$ in \eqref{fifth} gives us
    \[
        \lim_{j \to \infty} \int_{\Omega \cap f^{-1}(B_{r_j}(x_0))} u^{q-2} |\nabla\vec{H}|^2 \,d\mu = 0.
    \]
    After passing to an increasing sequence if necessary, these domains exhaust $\Omega$, which shows that
    \[
        \int_\Omega u^{q-2} |\nabla\vec{H}|^2 \,d\mu = 0.
    \]
    On $\Omega$, the weight $u^{q-2}$ is positive so $\nabla\vec{H} \equiv 0$ there. On the interior of $\{ u = 0 \}$, $\vec{H}$ is identically zero so $\nabla \vec{H} = 0$ there as well. At the boundary points of the zero set, the conclusion follows by the continuity of $\nabla \vec{H}$. Therefore, $\nabla\vec{H} \equiv 0$ on all of $M$.
\end{proof}

To finish up the proof of Theorem \ref{main}, we apply the following lemma in \cite{W13}:

\begin{lemma}[{see \cite{W13}, Lemma 3}] \label{glen lemma}
    Let $f: M^m \to \R^{m+n}$ be a smooth immersion whose mean curvature vector $\vec{H}$ satisfies $\nabla\vec{H} \equiv 0$, where $\nabla$ is the pull-back Euclidean connection. Then $\vec{H} \equiv 0$ and $f$ is minimal.
\end{lemma}

Combining Lemmas \ref{main lemma}, \ref{main lemma, weak} and \ref{glen lemma}, we obtain Theorem \ref{main}.

\section{A Liouville type result for $p$-biharmonic complete submanifolds of a Riemannian manifold}

In \cite{H15}, Han studied the $p$-biharmonic version of Chen's conjecture and proved that any complete \(p\)-biharmonic submanifold \((M,g)\) of a Riemannian manifold \((N,h)\) with non-positive sectional curvature satisfying  
\[
\int_M u^q \, d\mu < \infty
\qquad \text{for some } q\in(0,\infty)
\]
is minimal. In this section we consider $p$-biharmonic immersions for $p>1$. Recall that for an isometric immersion
\[
\varphi:(M,g)\longrightarrow (N,h),
\]
the necessary and sufficient condition for the $p$-biharmonicity of $\varphi$ is
\begin{equation}
-\widetilde{\Delta}
\bigl(u^{p-2}\vec H\bigr)
-\sum_{i=1}^m
R^N\bigl(u^{p-2}\vec H,e_i\bigr)e_i
=0.
\tag{2.1}
\end{equation}
Decomposing the above equation into its normal and tangential components, one can see that the $p$-biharmonicity is equivalent to
\begin{equation}\label{eh22}
    \Delta^\perp \bigl(u^{p-2}\vec H\bigr) - \sum_{i=1}^m A\left( e_i,S_{u^{p-2}\vec H}(e_i) \right) + \left( \sum_{i=1}^m R^N\bigl(u^{p-2}\vec H,e_i\bigr)e_i \right)^\perp =0,
\end{equation}
and
\begin{equation}\label{eh23}
    \operatorname{Tr}_g \left( \nabla^M S_{u^{p-2}\vec H} \right) + \operatorname{Tr}_g \left( S_{\nabla^\perp(u^{p-2}\vec H)}(\,\cdot\,) \right) - \left( \sum_{i=1}^m R^N\bigl(u^{p-2}\vec H,e_i\bigr)e_i \right)^\top =0.
\end{equation}
Here,
\[
    \Delta^\perp = \sum_{i=1}^m \left( \nabla^\perp_{e_i}\nabla^\perp_{e_i} - \nabla^\perp_{\nabla_{e_i}e_i} \right)
\]
is the Laplace operator associated with the normal connection $\nabla^\perp$. For $1<p<2$, the field $u^{p-2}\vec H$ need not be classically twice differentiable across $\{u=0\}$. We therefore use the equations above classically on $\Omega:=\{u>0\}$; in the proof of Theorem \ref{d1}, test functions vanishing near $\{u=0\}$ justify the passage from $\Omega$ to $M$. In \cite{H15}, Han also showed that if the sectional curvature is strictly negative, then the $L^q$-integrability of $u$ can be replaced by
\[
    \int_{B_r(x_0)} u^{a+2p-2}\, d\mu = O(r^d)
\]
for some $d > 0$. By utilizing the technique developed by Han, Cao and Luo \cite{CL16} later verified Chen's conjecture for $p$-biharmonic manifolds whose induced volume is of at most polynomial growth. In this section, we show that Han's result holds true for ambient manifolds whose sectional curvatures are nonpositive.
\begin{theorem} \label{d1}
Let $N$ be a Riemannian manifold whose sectional curvatures are nonpositive and let $f : M^m \to N^{m+n}$ be a complete $p$-biharmonic immersion for some $p > 1$. Suppose that, for some $q\geq0, d>0$,
\begin{equation} \label{new growth}
    \int_{B_r(x_0)} u^q\,d\mu\le Cr^d
\end{equation}
for sufficiently large $r$, where $B_r(x_0)$ is the intrinsic geodesic ball in $M$ centered at a fixed point $x_0$. Then $\vec{H}\equiv0$ and hence $f$ is minimal.
\end{theorem}

Before proceeding with the proof, we introduce a Bochner inequality for $u^{2p-2}$. Let $f:M^m\to N^{m + n}$ be an immersion, and denote by $\tau(f)$ its tension field. Then,
\[
\tau(f)=\operatorname{tr}A.
\]
with $A$ being the second fundamental form corresponding to $f$. In this paper, we use the convention
\[
\vec{H}=\operatorname{tr}A,
\]
which differs from that used in \cite{H15} by a factor of $\frac{1}{m}$. 

Recall that throughout this paper, we use the same convention for the curvature tensor as in \cite{H15}.

On $\Omega:=\{u>0\}$, set
\[
V = u^{p-2}\vec{H}.
\]
Since $f$ is $p$-biharmonic, $\tau_p(\varphi) = 0$. Hence, by \eqref{phan}, the normal component of the \(p\)-biharmonic equation becomes
\begin{equation}\label{e3.1}
\Delta^\perp V
=\sum_{i=1}^m A(S_Ve_i,e_i)
-\sum_{i=1}^m\bigl( R(V,e_i)e_i\bigr)^\perp.
\end{equation}
Here \(\{e_1,\ldots,e_m\}\) is a local orthonormal frame on \(M\), \(S_V\) denotes the shape operator associated with the normal vector field \(V\), and \(R\) is the curvature tensor of the ambient manifold.

On the other hand, the Bochner identity for the normal vector field $V$ gives
\[
\frac12\Delta |V|^2
=\left\langle\Delta^\perp V,V\right\rangle
+
|\nabla^\perp V|^2.
\]
This together with 
\eqref{e3.1} and the definition of $V$ yields
\begin{equation}\label{e3.4}
\frac12\Delta u^{2p-2}
=
|\nabla^\perp V|^2
+
\sum_{i=1}^m
\left\langle A(S_V e_i,e_i),V\right\rangle\
-
\sum_{i=1}^m
\left\langle R(V,e_i)e_i,V\right\rangle.
    \end{equation}
Let $\nu := \frac{1}{u} \vec{H}$, we then have $V = u^{p-1}\nu$ and $ S_V=u^{p-1}S_\nu$. From there, the second term on the right-hand side of \eqref{e3.4} can be estimated as follows. First, we have
\begin{align}
\sum_{i=1}^m
\left\langle A(S_Ve_i,e_i),V\right\rangle
&=
u^{2p-2}
\sum_{i=1}^m
\left\langle A(S_\nu e_i,e_i),\nu\right\rangle\notag\\
&=
u^{2p-2}
\sum_{i=1}^m
\left\langle S_\nu^2e_i,e_i\right\rangle\
=
 u^{2p-2}|S_\nu|^2. \label{e3.5}
\end{align}
Since \(S_\nu\) is a self-adjoint endomorphism of \(TM\) and
\[
\operatorname{tr}S_\nu
=\sum_{i=1}^m
\left\langle A(e_i,e_i),\nu\right\rangle
=\left\langle \vec{H},\frac{\vec{H}}{u}\right\rangle
=u,
\]
the Cauchy-Schwarz inequality for the eigenvalues of \(S_\nu\) gives
\[
|S_\nu|^2
\geq
\frac{1}{m}\bigl(\operatorname{tr}S_\nu\bigr)^2
=\frac{1}{m}u^2.
\]
This together with \eqref{e3.5} implies that
\begin{equation}\label{e3.7}
\sum_{i=1}^m
\left\langle A(S_V e_i,e_i),V\right\rangle
\geq
\frac{1}{m}u^{2p}.
\end{equation}
Now since the ambient manifold has non-positive sectional curvature, by using \eqref{e3.4} and \eqref{e3.7}, we finally obtain
\begin{equation}\label{e3.9}
\frac12\Delta u^{2p-2}
\geq
|\nabla^\perp V|^2
+
\frac{1}{m} u^{2p}.
\end{equation}

Next, we want to show a Kato-type inequality for $V$. Since $V = u^{p-2}\vec{H}$,
we have
\[
\nabla^\perp_{e_i}V = u^{p-2}\nabla^\perp_{e_i}\vec{H} + (p-2)u^{p-4} \left\langle\nabla^\perp_{e_i}\vec{H},\vec{H}\right\rangle \vec{H}.
\]
Therefore,
\begin{align*}
|\nabla^\perp_{e_i}V|^2
=&
u^{2p-4}
\left|\nabla^\perp_{e_i}\vec{H}\right|^2\\
&+2(p-2)u^{2p-6}
\left\langle\nabla^\perp_{e_i}\vec{H},\vec{H}\right\rangle^2+(p-2)^2u^{2p-8}
\left\langle\nabla^\perp_{e_i}\vec{H},\vec{H}\right\rangle^2 u^2\\
=& u^{2p-4} \left|\nabla^\perp_{e_i}\vec{H}\right|^2+p(p-2) u^{2p-6} \left\langle\nabla^\perp_{e_i}\vec{H},\vec{H}\right\rangle^2.
\end{align*}
Summing over \(i=1,\ldots,m\), we obtain
\[
|\nabla^\perp V|^2
=
u^{2p-4}|\nabla^\perp \vec{H}|^2
+
p(p-2)u^{2p-6}
\sum_{i=1}^m
\left\langle\nabla^\perp_{e_i}\vec{H},\vec{H}\right\rangle^2.
\]
Since $\left\langle\nabla^\perp_{e_i}\vec{H},\vec{H}\right\rangle^2\leq|\nabla^\perp_{e_i}\vec{H}|^2|\vec{H}|^2$, we obtain the following Kato-type inequality:
$$
    |\nabla^\perp V|^2 \geq \min\{1,(p-1)^2\} u^{2p-4}|\nabla^\perp \vec{H}|^2.
$$
From there, by letting $c_p:=\min\{1,(p-1)^2\}>0$, the Bochner formula \eqref{e3.9} implies
\begin{equation}\label{e3.11}
\Delta u^{2p-2} \geq 2c_p u^{2p-4}|\nabla^\perp \vec{H}|^2 + \frac{2}{m}u^{2p},
\end{equation}
wherever $\vec{H} \neq 0$.

The Bochner-type formula and Kato-type inequality that we just presented are valid in any codimension since we only use the shape operator \(S_\nu\), where \(\nu=\frac{1}{u}\vec{H}\) is the unit normal vector field in the direction of the mean curvature vector.

\begin{proof}[Proof of Theorem \ref{d1}]
Let $a>0$ and set $b:=a+2p>4$; we shall choose $a$ explicitly below. For each $r>0$, choose $\lambda_r\in C_c^\infty(M)$, with $0\leq\lambda_r\leq1$, such that $\lambda_r\equiv1$ on $B_r(x_0)$, $\lambda_r\equiv0$ outside $B_{2r}(x_0)$, and $|\nabla^M\lambda_r|\leq c/r$.

Let $\zeta_\sigma$ be the nondecreasing zero-set cut-off introduced in the proof of Lemma \ref{main lemma, weak}, and put
\[
    \psi_\sigma:=\lambda_r^b u^a\zeta_\sigma(u)\in C_c^1(\Omega).
\]
We may therefore test \eqref{e3.11} on $\Omega$ with $\psi_\sigma$. In the resulting integration by parts, the terms in which the derivative falls on $u^a$ or $\zeta_\sigma(u)$ are nonpositive, because $\zeta_\sigma'\geq0$ and
\[
    \nabla^M u^{2p-2}=2(p-1)u^{2p-3}\nabla^M u
    \qquad\text{on }\Omega.
\]
Estimating the remaining cut-off term using $|\nabla^M u|\leq|\nabla^\perp\vec H|$ and Young's inequality gives
\begin{align*}
&c_p\int_\Omega \lambda_r^b u^{b-4}\zeta_\sigma(u)|\nabla^\perp\vec H|^2\,d\mu
+\frac{2}{m}\int_\Omega \lambda_r^b u^b\zeta_\sigma(u)\,d\mu\\
&\qquad\leq C_{p,b}\int_\Omega
\lambda_r^{b-2}u^{b-2}\zeta_\sigma(u)|\nabla^M\lambda_r|^2\,d\mu.
\end{align*}
As $\sigma\downarrow0$, $\zeta_\sigma(u)\uparrow1$ on $\Omega$, Monotone Convergence gives
\begin{align}
&c_p\int_M
 \lambda_r^b u^{b-4}|\nabla^\perp\vec H|^2\,d\mu
 +\frac{2}{m}\int_M\lambda_r^b u^b\,d\mu                                  \notag\\
&\qquad\leq
 C_{p,b}\int_M
 \lambda_r^{b-2}u^{b-2}|\nabla^M\lambda_r|^2\,d\mu.       \label{global-test}
\end{align}
Here $b>4$ ensures that all three weighted integrands extend continuously by zero across $\{u=0\}$. For any $s\in(0,b-2]$, applying Young's inequality for the right-hand side of \eqref{global-test} with exponents $\frac{b}{s}$ and $\frac{b}{b-s}$ gives
\begin{equation} \label{absorption}
    \int_M u^{b - 2} \lambda_r^{b - 2} | \nabla^M \lambda_r|^2 \, d\mu \leq \epsilon \int_M u^b \lambda_r^b \, d\mu + C_\epsilon \int_M u^{\frac{(b-2-s)b}{b-s}} \lambda_r^{\frac{(b-2-s)b}{b-2}} | \nabla^M \lambda_r|^{\frac{2b}{b - s}} \, d\mu. 
\end{equation}
By choosing $\epsilon := \frac{1}{mC_{p,b}}$, we are able to absorb the first term on the right-hand side of \eqref{absorption} to the second term on the left-hand side of \eqref{global-test}. Then by using the gradient bound $|\nabla^M \lambda_r| < \frac{c}{r}$, 
\begin{equation} \label{d3}
    \begin{aligned}
        c_p\int_M \lambda_r^b u^{b-4} |\nabla^\perp \vec{H}|^2 \, d\mu +\frac{1}{m}\int_M \lambda_r^b u^b \, d\mu& \\
        \leq \left(\frac{C}{r}\right)^{\frac{2b}{b-s}} \int_{B_{2r}(x_0)} u^{(b-2-s)\frac{b}{b-s}} \, d\mu,& 
    \end{aligned}
\end{equation}
where $C$ is a positive constant depending on $a, p,$ and $s$.

Now fix $q\geq0$ and choose $a$ (and hence $b$) sufficiently large that
\[
    q < a + 2p - 2.
\]
Thus $q<b-2$. We claim that there exists $s \in (0 , b - 2]$ such that
\[
    q = \frac{b(b - 2 - s)}{b - s} = (a + 2p - 2 - s) \frac{a + 2p}{a + 2p - s}.
\]
Indeed, let
\[
  F(s)=\frac{b(b-2-s)}{b-s}.
\]
This function can be extended continuously to $[0,b-2]$, with
\[
  F(0)=b-2,
  \qquad
  F(b-2)=0,
  \qquad
  F'(s)=-\frac{2b}{(b-s)^2}<0.
\]
Hence, $F$ decreases strictly and maps $[0,b-2]$ to $[0,b-2]$. If $0 < q < b - 2$, the corresponding $s$ lies in $(0 , b - 2)$, while for $q = 0$, one has $s = b - 2$.

Substituting this choice of $s$ into \eqref{d3} yields
\begin{equation} \label{d2}
    \begin{aligned}
        \int_M \lambda_r^{b}u^{b-4} |\nabla^\perp \vec{H}|^2 \, d\mu +\frac{1}{m}\int_M \lambda_r^{b}u^{b}\, d\mu &\le \frac{C}{r^{b-q}}
        \int_{B_{2r}(x_0)} u^q \, d\mu \\
        &\leq \frac{Cr^d}{r^{a+2p-q}}.
    \end{aligned}
\end{equation}
For each fixed $R>0$, choose $r>R$. Since $\lambda_r\equiv1$ on $B_r(x_0)$, we also have $\lambda_r \equiv 1$ on $B_R(x_0)$. Therefore, \eqref{d2} gives

\[
  \frac1m\int_{B_R(x_0)}u^b\, d\mu
  \leq
  \frac1m\int_M\lambda_r^bu^b\, d\mu
  \leq
  C r^{d-(b-q)}.
\]
By choosing
\[
  a=q+d+4, \qquad b=a+2p=q+d+2p+4,
\]
we have
\[
  d-(b-q)=d-(d+2p+4)<0.
\]
Letting $r\to\infty$ while keeping $R$ fixed gives
\[
  \int_{B_R(x_0)}u^b\, d\mu=0.
\]
Finally, as $R \to \infty$ (passing to an increasing sequence, if necessary), the family of balls $B_R(x_0)$ exhausts $M$. This shows that
\[
    \int_M u^b \, d\mu = 0.
\]
Since the integrand $u^b$ is nonnegative, we conclude that $u = 0$. Therefore, by continuity, $\vec{H} \equiv 0$  everywhere on $M$ and $f$ is minimal.
\end{proof}
\begin{remark}
    Theorem \ref{d1} extends Han's polynomial-growth conclusion from the strictly negative curvature case to the nonpositive curvature case, under the growth assumption \eqref{new growth}.
\end{remark}
\begin{remark}
   For \(q=0\), the growth condition \eqref{new growth} is interpreted as the polynomial volume growth
   \[
        \operatorname{Vol}(B_r(x_0))\leq Cr^d,
   \]
   which was studied by Cao and Luo (see \cite{CL16}, Theorem 1.2). Thus, Theorem \ref{d1} may be seen as a generalization of their result.
\end{remark}

\section*{Acknowledgement}

The research of the second author is supported by NAFOSTED under grant number: 101.02-2025.62. The research of the first author is supported by Australian Research Council grant FT250100880, which provides the principal funding for his scholarship and tuition waiver, with additional support from the University of Wollongong. The research of the third author has received support from Australian Research Council grants FT250100880 and DP250101080. This work began during a visit by the third author to the Vietnam Institute for Advanced Study in Mathematics (VIASM), Hanoi; he thanks VIASM for its hospitality.

\end{document}